\documentclass[12pt]{amsart}

\usepackage{common}
\usepackage{mathrsfs}

\title{Minimal stable types in  Banach spaces}

\author{Saharon Shelah}
\author{Alexander Usvyatsov}

\address{Saharon Shelah\\
Mathematics Department\\
Hebrew University of Jerusalem\\
91904 Givat Ram, Israel\\
}
\address{ Saharon Shelah \\ Department of Mathematics\\
Hill Center-Busch Campus Rutgers, The State University of New Jersey
110 Frelinghuysen Rd Piscataway, NJ 08854-8019, USA}

\urladdr{http://shelah.logic.at}

\address{Alexander Usvyatsov\\
    Universidade de Lisboa \\
  Centro de Matem\'{a}tica e Aplica\cc{c}\~{o}es Fundamentais\\
  Av. Prof. Gama Pinto,2\\
  1649-003 Lisboa \\
  Portugal}

\address{Alexander Usvyatsov\\
Mathematics Department\\
Hebrew University of Jerusalem\\
91904 Givat Ram, Israel\\
}

\thanks{Publication no. 1020 on Shelah's list of publications.
    Shelah thanks the Israel Science Foundation
   (Grant 1053/11), and the
   European Research Council (Grant 338821)
   for partial support of this research.
}
\thanks{
    Research partially supported by Marie Sklodowska Curie CIG 321915 "ModStabBan". Usvyatsov thanks the FCT (Funda\cc{c}\~{a}o para a Ci\^{e}ncia e a Tecnologia) for partial support
    of this research:
    FCT grant no. SFRH / BPD / 34893 / 2007; FCT Research Project PTDC/MAT/101740/2008; FCT Research Project PTDC/MAT/122844/2010; Development grant IF/01726/2012. }

\date{\today}

\begin{document}

\begin{abstract}
    We prove existence of wide types in a continuous theory expanding a
    Banach space, and density of minimal wide types among stable types in such
    a theory. We show that every minimal wide stable type is
    ``generically''  isometric to an
    $\ell_2$ space. We conclude with a proof of the following formulation of Henson's Conjecture: every
    model of an uncountably categorical theory expanding a Banach space is prime over a spreading model, isometric to the standard basis of a
    Hilbert space.
\end{abstract}

\maketitle

\section{Introduction}




The main motivation for this work is a conjecture formulated by
C. Ward Henson in the 1970s concerning geometric structure of
non-separably categorical elementary classes of Banach spaces. Several years ago, after
some partial progress had been made on Henson's question, the second
author suggested a concrete formulation of the conjecture. In this paper, we prove a more general result. 
We establish a structure theorem for non-separably categorical Banach {structures}, 
that is, any continuous structure expanding an underlying Banach space. 
In addition, our techniques suggest the beginning of geometric structure theory for a
larger class of \emph{stable} elementary classes of Banach spaces.

Essentially, a (complete) metric structure $B$ of density character $\lam$ is called \emph{categorical} (in $\lam$) if any 
(complete) structure $B'$ of density $\lam$ that is elementary equivalent to $B$ is in fact isometric to $B$. So 
categoricity means that the isomorphism type of $B$ is captured in a strong way by ``first order axioms that $B$ satisfies'', or, more precisely, by the (continuous first order) theory of $B$. An alternative way of defining categoricity (that does not explicitly involve logic) is the following. $B$ as above is categorical in $\lam$ if: whenever $B'$ has the same density character as $B$, and $B$ and $B'$ have isometric ultra-powers, then $B$ and $B'$ are already isometric. In other words, in some sense, the isomorphism type of $B$ is essentially determined by its local structure (and its density character). 

One may wonder why we consider \emph{non-separable} categoricity, that is, categoricity in an uncountable density $\lam$. 
This has to do with deep model theoretic phenomena and the history of development of classical model theory, some of which we try to explain later in this section. A short answer is that separable (or countable, in the particular case of classical - that is, discrete - first order theories) categoricity arises for very different reasons, and does not lead to a similar structure theory. This phenomenon that may seem peculiar at first, is a particular case of a general principle in model theory, according to which, even if one is only interested in ``small'' objects, it is instructive and helpful to consider larger structures first. It turns out that it is \emph{non-separable} categoricity that captures the property of ``isomorphism type of a structure is determined by the axioms that is satisfies'' (as much as this is possible in the context of first order axioms that are preserved under taking ultra-products). It leads to a strong and beautiful structure theory, implying that the isomorphism type of our structure $B$ is completely determined by a certain (abstract) dimension. In this article we show that if $B$ is a continuous expansion of a Banach space (for example, a real or a complex Banach space, a Banach lattice, a Banach algebra, and so on), the underlying dimension is, in fact, quite concrete: it is the linear dimension of a naturally occurring Hilbert space, which essentially determines the structure of $B$. 



\medskip

We now phrase the main problem that motivated this work in a more precise way. Since we study categoricity in power, and would like to consider all structures elementary equivalent to a given structure $B$ simultaneously, it is more convenient speak in terms of elementary classes. Since the original (Henson's) problem was stated specifically for Banach spaces (and since we hope to make our presentation understandable to a general audience), we will for the moment restrict our attention to this 
more concrete context. However, as mentioned above, our results encompass a much wider spectrum of structures; see Subsection \ref{context}.   

Let us remind the reader that a class $K$ of Banach spaces is called
\emph{elementary} if it is closed ``nicely'' under the ultra-product construction. More precisely, $K$ is elementary if it is closed under isometries, ultra-products
and ultra-roots (the last condition means that 
the complement of $K$ is closed under ultra-powers). It is well known that analyzing ultra-products and ultra-powers of a Banach space can be helpful (and often essential) for understanding its local structure. This suggests that it is natural to consider a Banach space \emph{together with all its ultra-powers} -- that is, even if one is only interested in the 
geometry of a particular space, it can be instructive to look at the elementary class that it generates. Hence elementary classes of Banach spaces are objects of interest. 

Equivalently, a class of Banach spaces is
elementary if it can be axiomatized in an appropriate logic. One can
work with either Henson's logic of positive bounded formulae
\cite{Hen} or continuous first order logic \cite{BU,BBHU}. 

Many
``natural'' classes of Banach spaces are elementary, for example:

\begin{itemize}
\item Fix $1 \le p <\infty$. Then the class of all Banach spaces
  isometric to $L_P(\mu)$ for some measure $\mu$ is elementary.
\item The class of all Banach spaces whose dual is isometric to
  $L_1(\mu)$ for some measure $\mu$ is elementary.
\item The class of all Banach spaces isometric to $C(K)$ for some
  compact Hausdorff space $K$ is elementary. In this case the precise
  axiomatization is not known (but it has been shown that this class
  is closed under ultra-products and ultra-roots).
\end{itemize}

An elementary class of Banach spaces $K$ is called
\emph{categorical} in a cardinal $\lam$ if there is a unique $B \in
K$ of density character $\lam$ up to isometry. A class $K$ is called
\emph{uncountably categorical}, or \emph{non-separably categorical}, if it is categorical in some
uncountable $\lam$. The most basic example is the class of all
Banach spaces isometric to a Hilbert space. There are other known
examples, but in all of them the behavior of the class is
``controlled'' in a very strong sense by an underlying Hilbert
space.  

%

This led C. Ward Henson to make the following conjecture.

%


\begin{cnj}(Henson)
  Let $K$ be an uncountably categorical elementary class of Banach
  spaces. Then
  \begin{itemize}
  \item $K$ is categorical in all uncountable cardinalities.
  \item Any $B \in K$ of uncountable density character is ``very close'' to (and ``determined by'') an underlying Hilbert space.
  \end{itemize}
\end{cnj}

The first part of the Conjecture is simply an analogue of a well-known \L o\'{s}'s Conjecture (aka Morley's Theorem) in classical logic. It was established independently by the authors \cite{ShUs837}, and by Ita\"{i} Ben Yaacov \cite{Ben:morley}. Both proofs resembled classical proofs of analogous results in the first order context. In this paper we prove a version of the more interesting (the second) part of Henson's Conjecture.
Our main theorem is

\begin{thm}\label{main_theorem}
    Let $K$ be an uncountably categorical elementary class of Banach
  spaces (or, more generally, Banach structures).
  Then there is a separable $B_0 \in K$
  and a definable minimal wide type $p_0$ over $B_0$, such that 
	\begin{itemize}
		\item Any Morley sequence in
		  $p_0$  is isometric to the standard basis of a Hilbert space.
		\item Any non-separable $B \in K$ is
		  prime over a Morley sequence in $p_0$ (which is the fundamental sequence of a spreading model of $B_0$).
	\end{itemize}
%
\end{thm}

We explain the terms that appear in the statement later.

\smallskip

There are various improvements that one can make in the statement. For example, $B_0$ can be taken to be the countable saturated model of $K$.


\smallskip

We prove the theorem
above for any elementary class of Banach {structures}, that is, Banach spaces expanded with
continuous extra-structure (we explain the different contexts we
work in in section 2). However, we are currently unaware of any interesting natural examples of non-separably categorical non-trivial continuous expansions of a Banach space.

\smallskip

As a matter of fact, it is surprisingly hard to find categorical examples of Banach spaces, as well as to prove categoricity. Until quite recently, the only construction that was ``known to the experts'' to yield an uncountably categorical Banach space (the first written proof of this fact has appeared recently in Henson and Raynaud in \cite{HeRa}) was the following:

\begin{exm}\label{exm:finitedim}
    Let $E$ any finite-dimensional Banach space, and let $H(E)$ be the
    class of all Banach spaces which are isometric to the direct Hilbert sum of 
    an infinite dimensional Hilbert space with $E$.
Then $H(E)$ is elementary and categorical in all
    infinite density characters.
\end{exm}

Motivated by Theorem \ref{main_theorem}, Henson and Raynaud \cite{HeRa} have embarked on a journey of searching for new uncountably categorical Banach spaces. They have developed new techniques of proving categoricity in this context, which led them to discover many natural examples, more sophisticated than Example \ref{exm:finitedim}. More specifically, Henson and Raynaud have discovered a criterion that ensures that all models of a particular continuous theory of a Banach spaces are of the form $E \oplus_m H$, where $E$ is a separable (modular) base space, $H$ is a Hilbert space, and $\oplus_m$ is a modular direct sum. Under natural assumptions, this construction yields a non-separably categorical elementary class.

Many of the examples that come out of the work of Henson and Raynaud have $\aleph_0$ separable models, and exhibit a natural notion of dimension, finite or infinite; indeed, the existence of models of finite dimension gives rise to infinitely many separable models. It is natural to conjecture that the geometric analysis int this paper can be pushed further to cover the finite-dimensional case, hence recovering the appropriate analogue of the Baldwin-Lachlan Theorem, This, however, requires more work, and will be dealt with elsewhere. 

\smallskip

The following corollary, of Theorem \ref{main_theorem_GA} which is more accessible to the general audience, can be stated. 

\begin{thm}\label{main_theorem_GA}
  Let $K$ be an uncountably categorical elementary class of Banach
  spaces. Then any non-separable $B \in K$ is prime over a sequence isometric to the standard 
	basis of a Hilbert space (which is a spreading model of a fixed separable $B_0$ in $K$). 
\end{thm}

Let us point out that Theorem \ref{main_theorem_GA} is significantly easier than Theorem \ref{main_theorem}: an attentive reader should be able to deduce it directly from Dvoretzky-Milman Theorem (Fact \ref{Dvoretzky}) and compactness (the ``spreading model'' part requires a bit more, but is nevertheless straightforward). One of the main features of Theorem \ref{main_theorem} is that it gives a definable geometric object (a definable wide type) that generates the basis for the  Hilbert space that underlies any non-separable (in fact, any ``large enough'') member of $K$. One may ask whether stronger definability requirements hold: for example, whether the Hilbert space itself may be assumed to be a type-definable set. The answer is ``yes'' in all the examples that have been constructed so far. Another interesting question is whether $B_0$ can be chosen 
to be the prime model in $K$, in which case \emph{every} member of $K$ would be of the above form (a natural next step in the analysis of separable members in $K$).These seem to be natural directions for further research. 

\medskip

Let us explain some of the basic notions that appear in the statements above. 

A model $B \in K$ is called \emph{prime} over a set $A$ if whenever $A$ embeds into $B' \in K$ via $f\colon A \hookrightarrow B'$, there is an embedding of $B$ into $B'$ that extends $f$. A model $B_0 \in K$ is called \emph{prime} if it is prime over the empty set. Note that all the embedding in this case are \emph{isometries}. 

This notion may seem somewhat abstract, and it may not be clear why we believe that Theorem \ref{main_theorem} suggests that 
any (non-separable) $B \in K$ is \emph{determined} by the underlying Hilbert space $H$, given by the spreading model . It is therefore worth saying that one can make the relationship between $H$ and $B$ more concrete. that in this context $B$ is in fact ``constructible'' over $H$ in a certain sense. In fact, any element in $B$ realizes an \emph{isolated} type over $H$, which is a type that has to be realized in any model (element of $K$) that contains $B$. However, the precise meaning of these notions is somewhat technical, and will not be explained here. 

%

\smallskip

We also remind the reader that a sequence $\lseq{e}{i}{\lam}$ is called a \emph{spreading model} of a Banach space $B$ \cite{BrSu} if it is 1-subsymmetric (quantifier-free indiscernible), and there is a sequence $\lseq{b}{n}{\om}$
in $B$ which is asymptotically isometric to $\lseq{e}{i}{\lam}$: there exists a null sequence of positive reals $\lseq{\eps}{\ell}{\om}$ such that whenever $k>\ell$, we have 

$$ \left| \|\sum_{j=0}^{k-1}r_jb_{n_j} \| -\|\sum_{j=0}^{k-1}r_je_j\|\right|<\eps_\ell$$

for every $\ell<n_0<n_1<\ldots<n_{k-1}<\om$ and $r_j \in [-1,1]$.

Clearly, since \seq{e_n} is 1-subsymmetric, the sum $\sum_{j=0}^{k-1}r_je_j$ can be replaced with $\sum_{j=0}^{k-1}r_je_{i_j}$ for every $i_0<i_1<\ldots<i_{k-1}<\lam$.

\medskip

Equivalently, in model theoretic terms, a sequence $\lseq{e}{i}{\lam}$ is called a spreading model of $B$ if it is a quantifier free ``co-heir'' sequence over $B$: that is, it is quantifier free indiscernible, and $\tp(e_i/Be_{<i})$ is finitely satisfiable in $B$ (here ``$\tp$'' stand for the quantifier free type in the pure language of normed spaces).

\medskip

So Theorem \ref{main_theorem} shows that every non-separable $B \in K$ is essentially prime over a Hilbert space, which is nicely based on a ``small'' separable ``base'' space. 
In Example \ref{exm:finitedim} the base space is the direct sum of $E$ with a separable Hilbert space. 

Theorem \ref{main_theorem} is proven in section \ref{conclusions}, Theorem
\ref{thm:main}. 

\smallskip



Although this was not clear to us until the proofs were basically
finalized, a posteriori it has become apparent that this
formulation of Henson's Conjecture, and the techniques developed on
the way to its proof, provide in a sense a true analogue of
geometric characterizations of uncountably categorical elementary
classes in classical model theory, continuing the work of Baldwin,
Lachlan, Zilber and others (which we discuss in the next subsection)
in the context of Banach spaces. We believe that this paper lays the
foundations for the developing of geometric stability in this
setting. Hence it is our hope that the results here are not of
isolated interest, but rather a beginning of a new chapter in model
theoretic study of Banach spaces.

\subsection*{History and background}

Henson's Conjecture is strongly related to well-known results on
classical uncountably categorical elementary classes. In 1962 Morley
\cite{Mor} proved the conjecture of \L o\'s which stated that a countable first
order theory $T$ which is categorical in some uncountable power, is
categorical in any uncountable power. Basic examples of such
theories are the theory of algebraically closed fields of a fixed
characteristic, and the theory of vector spaces over a fixed countable field.
Morley's proof showed that an uncountably categorical theory $T$
admits a \emph{notion of independence} and that any model of $T$ is
both saturated (``rich'') and prime (``small'') over a basis with
respect to this notion.

Less than ten years later Baldwin and Lachlan \cite{BaLa} gave a different,
more geometric proof of Morley's Theorem. They showed that every
model of an uncountably categorical theory $T$ is determined by a
``strongly minimal'' definable set, on which the independence notion
is of a very special kind: it is determined by algebraic closure.
Their proof also gave information about countable models of
uncountably categorical  theories.

The results of Baldwin and Lachlan led to further research.
Specifically, Zilber studied geometric structure of strongly minimal
sets and showed that in many cases they are either ``field-like'' or
``group-like'' (and in the ``field-like'' case one can interpret an
algebraically closed field in the model). One reference for Zilber's
work is \cite{Zil}. A posteriori it turns out that Henson's
Conjecture called for a similar analysis for Banach spaces
(``interpreting'' a Hilbert space inside the model), but no
appropriate tools were available until very recently. For example,
no analogue of a strongly minimal set was known. In this article, we
introduce new geometric objects, which we call \emph{wide types}.
Our thesis is that \emph{minimal wide types} are an appropriate
analogue of strongly minimal sets in this setting.

Another important notion that we are going to make use of is
\emph{stability}. In his proof, Morley introduced the notion of
\emph{$\om$-stability}. He proved that an uncountably categorical
theory is $\om$-stable, and that $\om$-stability implies several good
properties, such as existence of prime models over any set (we shall
explain the notion of a prime model later) and a ``nice'' notion of independence. Later the first author defined
the more general notion of \emph{stability} and showed that
any stable theory admits a similar notion of independence.

Stability was first introduced to functional analysis by Krivine and
Maurey, who proved in \cite{KM} that any stable Banach spaces
contains an almost isometric copy of $\ell_p$ for some $p$. It was
further investigated by Iovino in \cite{Iov:stabI, Iov:stabII} and
other works, and, more recently, by Ben Yaacov and the authors (e.g.
\cite{Ben:morley, BU, ShUs837}).

We have already pointed out that the first part of Henson's Conjecture states
that the analogue of Morley's Theorem holds for classes of Banach
spaces. This was proven independently by the authors \cite{ShUs837}
and Ita\"{i} Ben Yaacov \cite{Ben:morley}. The two proofs are quite
different, but none of them gives much geometric information. In
some sense, both correspond to Morley's original proof, and do not
provide ``Baldwin-Lachlan analysis''. We will use several results
from \cite{ShUs837} in this article. In particular, we will use the
fact that uncountable categoricity implies a topological version of
$\om$-stability, which has property similar to those of classical
$\om$-stability. Consequently, uncountably categorical classes of
Banach spaces are stable.

We would also like to mention the classical theorem of Macintyre
\cite{Mac}: any $\om$-stable field is algebraically closed. In a
sense, this is a ``dual'' result to Morley Theorem: it shows
that ``algebraic'' structure follows from model-theoretic
properties. The second part of Henson's Conjecture has a similar flavor.

\subsection*{Acknowledgments}
We thank Ward Henson for numerous conversations that motivated and advanced this work. We are  grateful to Udi Hrushovski for many helpful comments,
and to Angus Macintyre for several inspiring conversations. We also thank the anonymous referee for helpful comments,
corrections, and suggestions.

\section{Preliminaries}

In this section we describe the framework in which we are going to
work. A reader who is familiar with continuous logic can easily skip
to the last subsection (subsection \ref{context}).

We refer the reader to \cite{Hei}, \cite{HI} or \cite{BBHU} for the
definition of an ultra-product of Banach spaces, and, more general,
normed structures.

In subsection \ref{basic case} we describe the very basic framework
of quantifier free formulas in the pure language of Banach spaces,
which is enough for proving Theorem \ref{main_theorem}. 
The
presentation in subsection \ref{basic case} should be accessible to
any mathematician, and quite familiar to a functional analyst. For
example, quantifier free types in this basic frameworks are
precisely what Krivine and Maurey defined as ``types'' in \cite{KM}.
Readers who desire to limit exposure to logic, can safely skip to section 3
directly after subsection \ref{basic case}.

In subsection \ref{general case} we present the more general context
of continuous logic expanding the Banach space structure. Working in
this framework, we prove more general results.

\bigskip
\subsection{The basic case}\label{basic case}

\begin{dfn}
  \begin{itemize}
  \item
    A \emph{quantifier free formula in the pure language of Banach
    spaces} over a set $A$ is an expression of the form $\|x+a\|$
    where $x$ is a variable and $a \in A$. We call such a formula
    \emph{a pure q.f.} formula.
  \item
    A \emph{pure q.f. condition} over $A$ is an expression of the form
    $\ph(x) = r$ where $\ph(x)$ is a pure q.f. formula over $A$ and
    $r \in \setR$.
  \item
    Let $\Sigma$ be a collection of pure q.f. conditions over a set
    $A$, $A \subseteq M$, $M$ a Banach space. We say that $\Sigma$ is
    \emph{approximately finitely satisfiable} in $M$ if for every
    finite $\Sigma_0 \subseteq \Sigma$ and $\varepsilon > 0$, there is
    $b \in M$ such that for every $[\ph(x) = r] \in \Sigma_0$ we have
    $\ph(b) \in [r-\varepsilon,r+\varepsilon]$.
  \item
    Given a Banach space $M$ and a subset $A \subseteq M$,
    a \emph{pure q.f. partial type} in $M$ over $A$ is a collection
    $\pi(x)$ of
    pure q.f. conditions over $A$ which is approximately finitely satisfiable
    in $M$, such that in addition $[\|x\|=r] \in \pi$ for some $r \in \setR$.
  \item Given a partial type $\pi(x)$, we say that the value of the formula
    $\ph(x)$ is \emph{determined} by $\pi$ if $[\ph(x) = r]\in \pi$
    for some $r \in \setR$. Otherwise we say that the value of $\ph$
    is \emph{undetermined} by $\pi$.
  \item Given a Banach space $M$ and a subset $A \subseteq M$, a
    \emph{complete pure q.f. type} in $M$ over $A$ is a partial pure
    q.f. type in $M$ over $A$ which determines the value of any pure
    q.f. free formula over $A$.

    In other words, a complete pure
    q.f. type $p$ over $A$ can be (and often is) viewed  as a function
    $\tau \colon A \to \setR$ such that for any $a \in A$ we have
    $\tau(a) = r$ if and only if $[\|x+a\|=r] \in p$.
  \item We denote the space of all complete pure q.f. types in $M$
    over $A$ by $\tS_{qf}(A,M)$, or just $\tS_{qf}(A)$ when $M$ is
    clear from the context.
  \item Given a partial type $\pi(x)$ and a formula $\ph(x)$, we
    denote by $\ph^\pi$ the \emph{value of $\ph$ according to
      $\pi$}. In other words, $\ph^\pi = r$ iff
    $[\ph(x)=r]\in\pi(x)$. Of course, this only makes sense if $\pi$
    determines the value of $\ph$.
  \item We say that $b \in M$ \emph{realizes} a partial type $\pi(x)$
    if $\ph(b) = \ph^\pi$ for every formula $\ph(x)$
  (whose value is determined by $\pi$).
  \end{itemize}
\end{dfn}

One can define complete types in a slightly different (equivalent) way, perhaps more
familiar to a non-logician. Let $M$ be a Banach space, and let $\fU$ be an ultrafilter on $M$.
Note that for any $a\in M$, the ultrafilter $\fU$ ``determines'' a value for the 
(quantifier free) formula $\|x+a\|$, which equals $\lim_{\fU} \|x+a\| $. In a sense, this is the 
``most likely'' (according to $\fU$) value of $\|x+a\|$ when $x$ varies over $M$.

Define the (quantifier free) \emph{average of $\fU$ over $M$} to be the complete (quantifier free) 
type $p$ over $M$ such that for every formula $\ph(x,a)$, its value according to $p$ equals $\lim_{\fU}\ph^p(x,a)$.
In particular, we have 

$$\|x+a\|^p = \lim_{\fU}\|x+a\| $$ 

for all $a \in M$. 

It should be very easy to see that this indeed gives rise to a complete quantifier free type over $M$. We denote 
it by $\Av_{qf}(\fU,M)$. Conversely, 
any complete type arises in this way: given a complete q.f. type over $M$, there is an ultrafilter $\fU$ on $M$ such that $p = \Av_{qf}(\fU,M)$. 

The latter remark allows us to characterize spreading models of $M$ in yet another way. Indeed, we say that a sequence 
$\lseq{a}{i}{\lam}$ is a sequence in a (q.f.) type $p \in S_{qf}(M)$ \emph{based on} $M$ if there exists an ultrafilter $\fU$ on $M$ such that $p = \Av_{qf}(\fU,M)$, and $a_i = \Av_{qf}(\fU,M\cup\Span\set{a_j\colon j<i})$ for all $i<\lam$. It is 
quite easy to see that (fundamental sequences for) spreading models of $M$ are exactly sequences of the form above, that is, 
sequences based on $M$ (in some $p \in S_{qf}(M)$). Such sequences are also called quantifier free co-heir sequences over $M$. 
It is easy to verify that this definition is equivalent to the ones given in the introduction.

\bigskip

The following version of Compactness Theorem can be found in e.g.
\cite{HI}.

\begin{fct}
  Let $M$ be a Banach space, $\pi(x)$ a partial type in $M$ over
  $A$. Then there exists an ultra-power $\hat{M}$ of $M$ and $b \in
  \hat{M}$ such that $b$ realizes $\pi$.
\end{fct}

\begin{dfn}
  We call a Banach space \emph{qf-saturated} if for every $A$ of
  cardinality less than the density character of $M$ and every $p \in
  \tS_{qf}(A,M)$, $p$ is realized in $M$.
\end{dfn}

Given an elementary class of Banach spaces, we will assume the
following:

\emph{There exists a Banach space $\fC$, which is qf-saturated, and
whose cardinality is much bigger than all other cardinals discussed
in this paper, and all $M \in K$, which are of interest to us, are
subspaces of $\fC$}.

Such $\fC$ is called the \emph{monster model} of $K$. There are
slight set-theoretic assumptions which are involved in the existence
of monster models, but we will not be concerned with these issues
here. In fact, in the cases that we are interested in in this paper
(e.g. if $K$ is uncountably categorical, or just stable), no such
assumptions are necessary.

\subsection{The general case}\label{general case}

In this subsection we give a very quick overview of continuous logic
in the special case of normed structures. The reader is referred to
\cite{BU} or \cite{BBHU} for details.

Just like in classical logic, a continuous \emph{signature} consists
of constant symbols, function symbols and predicate symbols. There
is a special predicate symbol for the norm, $\|\cdot\|$. Each
function symbol and predicate symbol is equipped with its arity $k
\in \setN$ and its modulus of uniform continuity, which is a
continuous function $\delta$ from $\setR_+$ to $\setR_+$ with
$\delta(0) = 0$. We will always assume that the signature contains
the signature of a vector space over $\setQ$; that is, it contains a
constant symbol $0$, a 2-ary function for vector addition, and for
every $q \in \setQ$, a 1-ary function $\cdot_q(x)$ for
multiplication by $q$.

A continuous \emph{pre-structure} $M$ for a given signature is a
semi-normed space, in which all the constant symbols are interpreted
as elements, function symbols - as functions on the structure,
predicate symbols - as functions from the structure to $\setR$. More
precisely, if $P$ is a predicate symbol of arity $k$, then its
interpretation $P^M$ is a function $P^M\colon M^k \to \setR$.
Similarly, if $f$ is a function symbol of arity $k$, then its
interpretation $f^M$ is a function $f^M\colon M^k \to M$.

Moreover, we demand that the predicate $\|\cdot\|$ is interpreted as
a semi-norm on $M$ and all the predicates and functions are
uniformly continuous with respect to $\|\cdot\|$, respecting their
continuous moduli. This ensures that the predicates and functions
are continuous uniformly over all structures. Roughly speaking, this
is what is needed in order to make ultraproducts work.


A \emph{structure} is a pre-structure in which $\|\cdot\|$ is a
complete norm.

One notion which is important to understand in order to read the
paper in full generality is that of a \emph{formula}. The algebra of
formulas is obtained as follows. An \emph{atomic} formula is an
expression of the form $P(\tau_1, \ldots, \tau_k)$ where $P$ is a
predicate symbol or arity $k$, and every $\tau_i$ is a \emph{term},
which is a ``generalized'' function symbol (an expression that can
be obtained by composing existing function symbols and applying them
to variables and constants).

For example, quantifier free formulas discussed in subsection
\ref{basic case}, which are expressions of the form $\|x+y\|$,
$\|x+a\|$ (where $a$ is a constant) or, more generally,
$\|\sum_{i<n}q_ix_i + a\|$ (where $q_i \in \setQ$, $x_i$ are
variables) are atomic formulas. Note that $qx$ means $\cdot_{q}(x)$,
so we omit the formal function symbol and use the familiar notation.

Now the algebra of formulas is the closure of the collection of
atomic formulas under ``connectives'' - bounded continuous function
from $\setR^k \to \setR$ (for some $k \in \setN$), ``quantifiers''
$\sup_x$ and $\inf_x$ (where $x$ is a variable) and uniform limits.
Note that due to uniform limits we obtain formulas of the form
$\|rx\|$ where $r \in \setR$, and due to connectives we can for
example speak of a formula $r\cdot\|x\|$ or $\left| \|x\| -
r\right|$, where $r \in \setR$. Using quantifiers, we get formulas
of the form
\begin{gather*}
    \sup_x \left| \|x\| - r \right|
\end{gather*}
The collection of all formulas (for a given signature) is also
called a \emph{language}.

A (closed) \emph{condition} is an expression of the form $[\ph \in
C]$, where $\ph$ is a formula and $C$ is a closed subset of $\setR$.
We will only work with conditions where $C$ is a closed interval,
often a point (most of the time $C = \set{0}$).

A variable in a formula $\ph$ is called \emph{bounded} if it is in a
scope of a quantifier, and it is called \emph{free} if it is not
bounded. Given a formula $\ph$ with free variables $x_1, \ldots,
x_k$, we often write $\ph(x_1, \ldots, x_k)$ in order to emphasize
the free variables. It is easy to see that a formula $\ph(x_1,
\ldots, x_k)$ and a structure $M$, defines a function $\ph^M \colon
M^k \to M$. In fact, $\ph^M$ is uniformly continuous, and, moreover,
uniformly so in all structures (one can calculate the uniform
continuity modulus of $\ph$, given the moduli of all function and
predicate symbols in the signature). So given a formula $\ph(x_1,
\ldots, x_k)$, a structure $M$, and $a_1, \ldots, a_k \in M$, one
can calculate $\ph^M(a_1, \ldots, a_k) \in \setR$. Hence given a
condition $[\ph(x_1, \ldots, x_k) \in C]$, a structure $M$, and
$a_1,\ldots,a_k \in M$, it makes to ask whether the condition is
\emph{true} in $M$ (denoted by $M \models [\ph \in C]$). If $M
\models [\ph \in C]$, we also say that $M$ is a \emph{model} of
(for) this condition.

A \emph{theory} is a collection of conditions with no free
variables, which has a model. We normally assume that a theory $T$
is closed under entailment, that is, if a condition $[\ph \in C]$
follows from $T$ (which means that it is true in all models of $T$),
then $[\ph \in C] \in T$. Compactness Theorem (see \cite{HI, BBHU})
states that a collection of conditions has a model if and only if
every finite subset of it does. A theory is called \emph{complete}
if for every condition $[\ph \in C]$ either it is in $T$ or for some
closed $D \subset \setR$ disjoint to $C$ we have $[\ph\in D]\in T$.
Equivalently, $T$ is complete if it ``forces'' a value for every
formula $\ph$ with no free variables, that is, $[\ph=r] \in T$ for
some $r \in \setR$. We will denote that value by $\ph^T \in \setR$.
Note that every theory can be extended to a complete theory (in
fact, every model $M$ of $T$ determines a complete theory).

We will normally assume that we have a fixed complete theory in the
background, and all structures are models of $T$; we will therefore
often simply call them ``models''. Given a model $M$, and a subset
$A$ of $M$, we will often expand the language by adding constant
symbols for all elements of $A$. Call this language $L(A)$. Then $M$
naturally becomes an $L(A)$-structure; we will call $L(A)$-formulas
``formulas over $A$''.

The next definition is of central importance. A \emph{type} $\pi(x)$
in a model $M$ over a set $A$ is a collection of conditions of the
form $\ph(x) \in [r_\ph,s_\ph]$ (where $\ph$ is a formula over $A$,
$r_\ph,s_\ph \in \setR$), which is finitely approximately
satisfiable in $M$. The latter means that for every finite subset
$\pi_0(x)$ of $\pi(x)$ and for every $\eps>0$ there exists $a \in M$
such that $\ph^M(a) \in [r_\ph-\eps,s_\ph+\eps]$ for every condition
$\ph(x) \in [r_\ph,s_\ph]$ in $\pi_0(x)$. Equivalently, by
Compactness, a type $\pi(x)$ is a collection of conditions of the
form above such that there is an ultrapower $\hat M$ of $M$ and $a
\in \hat M$ which satisfies all the conditions in $\pi(x)$. We say
that $a$ realizes $\pi$ and write $a \models \pi$.

In general, $x$ in the definition of the type does not need
to be a singleton (so neither does $a$, that is, maybe $a \in M^k$
for some $k \in \setN$), although in this paper it we will normally
work with formulas and types in one variable.

We say that a type $\pi(x)$ \emph{determines} a value of a formula
$\ph(x)$ if $[\ph(x) = r] \in \pi$ for some $r \in \setR$. A
\emph{complete} type over $A$ is a type over $A$ which determines a
value for every formula over $A$ (with the right number of
variables). We will denote the value of a formula $\ph$ ``according
to the type $\pi$'' by $\ph^\pi$.

There is a correspondence between complete theories and complete types
(a complete type can be viewed as a complete theory in an expanded
language). We denote the space of all complete type over a set $A$
in $n$ variables by $\tS^n(A)$. This is a compact Hausdorff
topological space, but we will not be concerned with this fact here.
Let $\tS(A) = \cup_{n<\om}S^n(A)$. Note that the space of types is
defined relatively to a certain model which contains $A$; but as we'll see in
a bit, we will be working in one big model of the theory $T$ (the
``monster'' model), and all types will be computed in that
structure.

Given a model $M$, a set $A$ and a tuple $a \in M^k$, we denote by
$\tp(a/A)$ the collection of all closed conditions over $A$ that $a$
satisfies. It is easy to see that $\tp(a/A) \in \tS(A)$, and we call
it the type of $a$ over $A$ (again, we forget to mention $M$).
Conversely, every complete type over a set $A$ is the type of some
$a$ over $A$ (possibly $a$ is in some ultrapower of $M$; soon this
won't matter because in the ``monster'' model we will have
realizations for all types over ``small'' sets).

Given a cardinal $\lam$, a model $M$ is called $\lam$-saturated if
every type over a subset of $M$ of cardinality less than $\lam$ is
realized in $M$. A model is called \emph{saturated} if it is
$|M|$-saturated. There is a mild set-theoretic assumption that goes
into the existence of saturated models, and it can be avoided if one
works with a slightly weaker notion than saturation (which has the
same properties that we care about), but we will not go into the
details here. As a matter of fact, in the cases that we will be
interested in in this paper (e.g. $T$ uncountably categorical, or
just stable), saturated models provably exist. Given a (complete)
theory $T$, we will assume the following:

\emph{There exists a saturated model $\fC$ of cardinality $\ka^*$
for some big enough cardinal $\ka^*$, that is, much bigger than all
cardinals mentioned in this paper (except $\ka^*$ itself, of
course). We call $\fC$ the ``monster model'' of $T$}.

\medskip

A useful consequence of saturation is the following homogeneity
property of the monster model: given two tuples $a,b \in M^k$ and a
set $A$ (of ``small'' cardinality, that is, less than $\ka^*$),
$\tp(a/A) = \tp(b/B)$ if and only if there exists $\sigma \in
\Aut(\fC/A)$ (the group of automorphisms of $\fC$ fixing $A$
pointwise) such that $\sigma(a) = b$.

Another useful notion (although we will not really need it here) is
that of an \emph{elementary submodel}: if $M$ is a substructure of
$N$, we say that $M$ is \emph{elementary} in $N$, $M \prec N$, if
for any formula $\ph$ over $M$ with no free variables, we have
$\ph^M = \ph^N$. For example, $M$ is always elementary in any of its
ultrapowers (this is \L o\'s's Theorem adapted to this context; see
\cite{BU, BBHU}).

\bigskip

The monster model of $T$ embeds elementarily any $M \models T$ of
``small'' cardinality. This is why we will be able to assume that
all models of $T$ are elementary submodels of $\fC$. Moreover, types
over subsets of $M$ are the same in $M$ and any elementary
extension; so it will be enough to talk about types in $\fC$ (and we
will not mention it).

\subsection{$\Delta$-types}\label{sub:Deltatypes}

The following definitions and notations are somewhat less standard than what is mentioned in the previous subsection, and are used extensively throughout the paper. 

Let $\Delta$ be a collection of formulas (with no parameters). We say that $\ph(x,a)$ is a \emph{$\Delta$-formula} if $\ph(x,y) \in \Delta$. A (partial) $\Delta$ type (over a set $A$) is a partial type that consists entirely of $\Delta$-formulae. A \emph{complete $\Delta$-type} over$A$ is a $\Delta$-type that determines a value for every $\Delta$-formula over $A$. We denote the set of all complete $\Delta$-types over $A$ by $\tS_\Delta(A)$. Similarly, we define the $\Delta$-type of an element (tuple) $a$ over a set $A$; it will be denoted by $\tp_\Delta(a/A)$. Just like with ordinary complete types, a complete $\Delta$-type is a $\Delta$-type of an element (tuple), and vice versa.

Normally  $\Delta$ will be assumed to be closed under connectives and substitution of variables. We will call such subsets of formulas \emph{fragments} of the language. 

If the language of $T$ expands the language of Banach spaces, then the set $\Delta = \Delta_{pqf}$ of all quantifier free formulas in the language of Banach spaces is a fragment of the language. In this case $S_\Delta$ is essentially what we called $S_{qf}$ in subsection \ref{basic case}.

\subsection{Context}\label{context}

The general context: $T$ is a continuous theory, whose monster model $\fC$
expands a Banach space $\fB$. We denote the language of $T$ by
$L=L_\fC$ and the language of Banach spaces (which is a part of $L$)
by $L_\fB$.

As we have mentioned before, one can restrict oneself to the following
context: $K$ is an elementary class of Banach spaces, $\fC = \fB$ is
its monster model.

As usual, all sets and tuples mentioned in the paper are subsets of
$\fC$ (of cardinality less than $|\fC|$), and all models are
elementary submodels of $\fC$ (again, of ``small'' cardinality).

As we have mentioned above, all types are types in $\fC$, and we will not mention this.

\section{Wide types over Banach spaces}

Recall that $\fC$ expands a real Banach space $\fB$.

\begin{dfn}
    We call a partial type in 1 variable $\pi(x)$ (possibly with
    parameters) \emph{wide} if the set of realizations of $\pi(x)$
    in $\fC$ contains the unit sphere of an infinite dimensional subspace of $\fB$.
\end{dfn}

\begin{rmk}
    The type $x=x$ is wide.
\end{rmk}


%

The main goal of this section is showing that complete wide types
exist over any set. We will make use of the following well-known
result, which is sometimes referred to as Concentration of Measure Phenomenon, or the Dvoretzky-Milman-Ramsey
Phenomenon. 
It is a consequence of the renown Dvoretzky's Theorem \cite{Dv}, but the approach we take is due to
Milman, e.g. \cite{Mil}, Theorem 1.2. We will refer to this fact as
Dvoretzky-Milman Theorem.

\begin{dfn}
\begin{enumerate}
\item
    Let $B$ be a Banach space, $\mathbf S(B)$ the unit sphere of $B$, $f
    \colon \mathbf S(B) \to \setR$. The \emph{spectrum} $\gamma(f)$ is the
    collection of all $r \in \setR$ such that for every $\eps>0$ and any integer $k$ there exists a
    $k$-dimensional subspace $F$ of $B$ such that $|f(x) - r| \le
    \eps$ for all $x$ in the unit square of $F$.
\item
    Let $B$, $f$ be as before. We denote by $\ga'(f)$ the collection
    of all $r \in \setR$ such that for any $k$ and $\eps$ as above,
    $F$ can be chosen to be $(1+\eps)$-isomorphic to a
    $k$-dimensional Hilbert space.
\end{enumerate}
\end{dfn}

\begin{fct}\label{Dvoretzky}(Dvoretzky-Milman Theorem).
    Let $f$ be a uniformly continuous function on the unit sphere of an infinite
    dimensional Banach space $B$. Then the spectrum $\ga(f)$ is not
    empty. Moreover, $\ga'(f)$ is not empty.

\end{fct}
\begin{prf}
    For the proof we refer the reader to e.g. \cite{BL}, section 12
    (specifically, combine Theorem 12.10 and Proposition 12.3 there).
	Alternatively, see \cite{Pes} for a detailed discussion of concentration phenomena. 
\end{prf}

The first approximation to our goal is the following.

\begin{prp}\label{prp:wideextend}
    Let $\pi(x)$ be a wide partial type (over a set $A$),
    $\ph(x,\a)$ be a formula. Then there exists $r \in \setR$ such that the partial type $\pi(x) \cup [\ph(x,\a) = r]$
    is wide.
\end{prp}
\begin{prf}
    Without loss of generality we may assume that $\|x\|^\pi = 1$.

    Let $B$ be an infinite dimensional subspace of $\fB$ whose unit
    sphere $\mathbf S(B)$ is contained in $\pi^\fC$. The formula $\ph(x,\a)$
    induces a uniformly continuous function $f$ from $\mathbf S(B)$ to $\setR$.
    By Dvoretzky-Milman Theorem \ref{Dvoretzky},
    $\ga'(f) \neq \emptyset$. Let $r \in \ga'(f)$.

    Let $H = \ell_2$. For every $v \in H$, introduce a free variable
    $x_v$. Let $\bold x = \inseq{x}{v}{H}$ Denote by $\Lambda(\bold x)$
    the linear quantifier free diagram of $H$ with variables $x_v$. That
    is,
    $$\Lambda(\bold x) = \set{x_v = \sum_{i<k}\lam_ix_{v_i} \colon v,v_i \in H,
    \lam_i \in \setR, v = \sum_{i<k}\lam_iv_i}$$

    Consider the following collection of formulas. This is the
    (approximate) quantifier free diagram of $H$ with the additional
    requirement that the unit sphere $\mathbf S(\mathbf x)$ satisfies $\pi(x) \&
    [\ph(x,\a) = r]$.

    \begin{gather*}
       \Gamma(\bold x) = \Lambda(\bold x) \cup \set{\pi(x_v) \& |\ph(x_v,\a) - r| \le \eps \colon \|v\|_H =
        1, \eps>0} \cup \\
        \left\{(1-\eps)\|v\|_H \le \|x_v\| \le (1+\eps)\|v\|_H  \colon
        v \in H, \eps>0 \right\}
    \end{gather*}

    We claim that $\Gamma(\bold x)$ is finitely satisfiable in $\fC$.
    Indeed, in order to make sure this is true, one has to argue that
    for any $k$ and any $\eps>0$ there is a $k$-dimensional subspace $F$
    of $B$ which is $(1+\eps)$-isomorphic to the $k$-dimensional Hilbert
    space $\ell^k_2$ such that $\ph(x,\a) \sim_\eps r$ on $\mathbf S(F)$, and
    this follows immediately from the fact that $r \in \ga'(f)$.

    This shows that $\pi(x) \cup \set{\ph(x,\a) = r}$ is a wide type (in
    fact, it contains the unit ball of an infinite dimensional subspace
    isometric to $\ell_2$).
\end{prf}

Modifying the proof of the Proposition above, we also obtain the
following.

\begin{lem}\label{lem:wideunion}
    Let \lseq{\pi}{i}{\lam} be an increasing chain of wide partial types. Then
    $\pi = \bigcup_{i<\lam}\pi_i$ is wide.
\end{lem}
\begin{prf}
    We use compactness as in the proof of Proposition
    \ref{prp:wideextend}. That is, let $\Lambda(\bold x)$ be as
    there, and let
\begin{gather*}
    \Gamma(\bold x) = \Lambda(\bold x) \cup \set{\pi_i(x_v) \colon \|v\|_H =
    1, i<\lam} \cup \\
    \left\{(1-\eps)\|v\|_H \le \|x_v\| \le (1+\eps)\|v\|_H  \colon
    v \in H, \eps>0 \right\}
\end{gather*}

Clearly $\Gamma$ is finitely satisfiable, hence consistent, so the
union $\pi$ is wide.

\end{prf}

\begin{thm}\label{thm:wideexist} (Existence of Wide Types).
Let $\pi(x)$ be a wide partial type over a set $A$, $\Delta$ a collection of formulae closed under
connectives. Then there exists a complete wide $\Delta$-type $p$ over $A$
extending $\pi$.
\end{thm}

\begin{rmk}
\begin{enumerate}
\item
	Recall the notion of a $\Delta$-type from Subsection \ref{sub:Deltatypes}.
\item
    We will normally use $\Delta = L$ or $L_\fB$ or $\Delta =$ quantifier free
    formulae in $L$ or $L_\fB$.
\item
	In case $\Delta = \Delta_{pqf}$ (quantifier free formulae in $L_\fB$), as noted in Subsection \ref{sub:Deltatypes}, we are in 
	the ``basic context'' described in Subsection \ref{basic case}. Hence whenever a subset $\Delta$ of the language is mentioned, 
	the reader can safely assume $\Delta = \Delta_{pqf}$, and we are working in the basic context. In particular, $\Delta$-types (partial and complete) are simply quantifier free types described in Subsection \ref{basic case}, and $S_\Delta(A) = S_{qf}(A)$. 
\item
    Note that since $x=x$ is wide, the theorem implies in particular that there exists a
    complete wide type over any set.
\end{enumerate}
\end{rmk}

\begin{prf}
    Without loss of generality we may assume that $\|x\|^\pi = 1$.




Enumerate all $\Delta$-formulae over $A$
\seq{\ph_\al(x,\a_\al)\colon \al < \lam = |A|+|T|} such that
\begin{itemize}
\item[(*)]
    If $\delta$ is a limit ordinal and $\al_1 < \al_2 < \ldots \al_k
    < \delta$, then for any $k$-ary connective $F$, for some $\al<\delta$ we
    have
    $$F(\ph_{\al_1}(x,\a_{\al_1}),\ldots,\ph_{\al_k}(x,\a_{\al_k})) =
    \ph_{\al}(x,\a_{\al})$$
\end{itemize}

Now construct an increasing continuous sequence of wide types
$\pi_\al$ by induction on $\al$ such that:
\begin{itemize}
\item
    $\pi_0(x) = \pi(x)$
\item
    $\pi_\al(x)$ determines the value of $\ph_\be(x,\a_\be)$ for all
    $\be<\al$
\end{itemize}


For successor stages, use Proposition \ref{prp:wideextend}, and for
limit stages apply compactness as in the proof of Lemma
\ref{lem:wideunion}. This is possible by (*) above.

Obviously $p = \pi_\lam$ is as required.

\end{prf}

Analyzing the proof, we see that we have actually shown

\begin{cor}
    Let $\pi(x)$ be a wide partial type, $A$ a set containing the domain
    of $\pi$, $\Delta$ a collection of formulae closed under
    connectives. Then there exists a complete $\Delta$-type $p$ over $A$
    containing $\pi$ such that $\pi^\fC$ contains the unit sphere of an infinite dimensional subspace isometric to a Hilbert space.
\end{cor}

\section{Wide stable types}

%
%


Let $\Delta$ be a fragment of the language (see Subsection \ref{sub:Deltatypes}).

In this section we will use the notion of the algebraic closure of a set
$A$, $\acl(A)$, which is the collection of all $b$ whose orbit under the action 
of the automorphism group $\Aut(\fC/A)$ is compact. Recall that any model is algebraically
closed, that is, $\acl(M)=M$.

%

Recall that a complete $\Delta$-type $p$ is called \emph{definable}
if for every $\Delta$-formula $\ph(x,y)$, the function $\theta(y) =
d_p x \ph(x,y)$ defined as $d_p x \ph(x,a) = \ph^p(x,a)$ is a
definable predicate (that is, can be uniformly approximated
by formulae). We say that $\pi$ is $\Delta$-definable if for
every $\Delta$-formula $\ph$, the function $d_p x \ph(x,y)$ can be
uniformly approximated by $\Delta$-formulae.

As an example, let $\Delta$ be the collection of all quantifier free
formulae in the language $L_\fB$, and let $p$ be a complete
quantifier free 1-type over a closed subspace $A$. So $p$ is
determined by conditions of the form $\|x+a\| = r_a$ for all $a \in
A$. In other words, $p$ is determined by the function $\tau_p \colon
A \to \setR$ defined by $\tau(a) = \|x+a\|^p$. We call $p$ definable
if this function is a definable predicate (that is, can be uniformly
approximated by formulae), and we call it quantifier-free definable,
if it can be uniformly approximated by quantifier-free formulae. 

Definability is quite a strong assumption; we elaborate on its meaning a little bit in Remark \ref{rmk:uniqueextension}.


We will not use the notion of a stable formula (as defined in
\cite{BU}) in this article. Let us just remark that a formula
$\ph(x,y)$ is stable if and only if every every $\ph$-type is
$\Delta$-definable, where $\Delta$ is a the closure of $\ph$ under
connectives and permutations of variables. For example, the norm
$\|x+y\|$ is stable in $\fC$ (in the sense of Krivine and Maurey
\cite{KM}) if and only if every quantifier-free type is
quantifier-free definable. See \cite{BU} (or \cite{Iov:stabI} for a
slightly less general formulation).

The following definition is a straightforward generalization of the
classical concept due to Lascar and Poizat:

\begin{dfn}
    Let $\Delta$ be a collection of formulae closed under connectives and permutations of variables. A partial $\Delta$-type
    $\pi$ (possibly with parameters) is called \emph{$\Delta$-stable} (or simply
    \emph{stable} when $\Delta$ is clear from the context) if every
    extension of it  to a $\Delta$-type over $\fC$ is $\Delta$-definable.
\end{dfn}

\begin{rmk}
    In \cite{Iov:def} Jos\'{e} Iovino studies quantifier free types over
    Banach spaces that he calls ``stable''.
    We would like to alert the reader to the fact that Iovino's concept is significantly weaker than the
    classical notion defined above. In a dependent theory \cite{Sh:c, Sh715} (if one considers all
    formulae, and not just q.f. ones), Iovino's definition  is
    equivalent to a (much more general than stability) notion of \emph{generic
    stability} \cite{Sh715,Us, HP}. In an arbitrary theory Iovino's definition
    is even weaker than generic stability: e.g., a $c_0$-type is Iovino-stable, but not generically stable (for the discussion 
	of generic stability in the general context see e.g. \cite{PT, GOU}).
\end{rmk}

We would like now to define \emph{forking}. The following definition
is equivalent to the classical one when one restricts attention
to stable types.

\begin{dfn}
    Let $p\in \tS_\Delta(A)$ be a complete stable $\Delta$-type over $A$, and let $\rho$ be a
    partial $\Delta$-type extending $p$ (so $\rho$ is stable as well). We say that $\rho$ \emph{does not
    fork over $A$} or \emph{is a non-forking extension} of $p$ if
    $\rho$ is definable over $\acl(A)$.

    If $\rho$ is not definable over $\acl(A)$, we say that it forks
    over $A$ (or is a forking extension of $p$).
\end{dfn}


The following is a classical fact about stable types (a straightforward generalization of \cite{LP} to 
the continuous context).

\begin{fct}\label{stationary}
    A complete stable type over an algebraically closed set is
    stationary, which means that it has a unique non-forking
    extension to a complete type over $\fC$.
%
\end{fct}

\begin{fct}\label{local character}
    Let $p = p_0$ be a $\Delta$-stable type. Then there does not exist an
    increasing sequence of $\Delta$-types \lseq{p}{i}{|\Delta|^+} such that
    $p_{i+1}$ is a forking extension of $p_i$.
\end{fct}
\begin{prf}
    Denote $\lam = |\Delta|$. Let $q = \cup_{i<\lam^+}p_i$. Since $q$ extends $p = p_0$ and $p$ is stable, $q$
    is $\Delta$-definable (hence definable over a subset $B$ of
    $\dom(p)$ of cardinality $\lam$). Clearly $B \subseteq
    \dom(p_i)$ for some $i$; but since $p_{i+1} = q \rest
    \dom(p_{i+1})$, this implies that $p_{i+1}$ is definable over
    $B$, hence is a nonforking extension of $p_i$, a contradiction.
\end{prf}

\begin{dfn}
	We call a wide partial $\Delta$-type $\pi(x)$ over a set $B$ \emph{wide $\Delta$-minimal} if $\pi$ has a
    unique extension to a global wide $\Delta$-type (a complete wide $\Delta$-type over $\fC$).
\end{dfn}

\begin{rmk}
\begin{enumerate}
\item We omit $\Delta$ when it is clear from the context. 
\item Note that by Theorem \ref{thm:wideexist} any wide partial type has at least one global wide extension. 
\end{enumerate}
\end{rmk}

From now on, let us fix $\Delta$ containing the quantifier free
formulae of $L_\fB$, closed under connectives and permutations of
variables. When we say ``type'', ``formula'', etc, we mean
$\Delta$-type, $\Delta$-formula.

\begin{cor}\label{cor:minexist} (Density of minimal types)
    Let $\pi(x)$ be a partial wide type over a set $A$. Then there exists
    $B \supseteq A$ with $|B\setminus A| \le |\Delta|$ and $p \in
    \tS_\Delta(B)$ which extends $\pi$ and is {wide minimal}. Moreover, the unique wide extension to $p$ to a global $\Delta$-type is the unique non-forking extension of $p$. In other words, no forking
    extension of $p$ to a ($\Delta$-type over a) superset of $B$ is wide.
\end{cor}
\begin{prf}
    Construct by induction an increasing continuous sequence of sets $A_i$ and  an increasing sequence of types $p \in
    \tS_\Delta(A_i)$ such that
    \begin{itemize}
    \item
        $A_0 = A$
    \item
        $|A_{i+1}\setminus A_i|$ is finite
    \item
        $p_0$ extends $\pi$
    \item
        $p_i$ is wide for all $i$
    \item
        $p_{i+1}$ forks over $A_i$
    \end{itemize}
    Successor stages of the construction are clear. For limit stages, use
    Lemma \ref{lem:wideunion}.
    Since the construction has to get stuck at some $i<|\Delta|^+$,
    clearly (by stationarity)
    $B = \acl(A_i)$ and any extension of $p_i$ to $B$ are as required.
\end{prf}

We will now study the structure of minimal wide (stable) types.

Let $O$ be a linearly ordered set. Recall that a sequence $I =
\inseq{a}{i}{O}$ is called $\Delta$-\emph{indiscernible} over a set
$A$ if the $\Delta$-type of any finite sequence $a_{i_1}\ldots
a_{i_k}$ over $A$ depends only on the order between the indices
$i_1, \ldots, i_k \in O$. So if $\Delta$ is the collection of all
the quantifier free formulae in the language $L_\fB$ and $A =
\emptyset$, then $I$ is $\Delta$-indiscernible if and only $I$ is
1-subsymmetric. As mentioned before, we will omit $\Delta$.

A sequence $I$ as above is called an \emph{indiscernible set} over
$A$ if the type of any finite sequence $a_{i_1}\ldots a_{i_k}$ over
$A$ depends only on the number $k$. As an example, one may think of
the standard basis of $\ell_p$.

The following is another classical fact about stable types:

\begin{fct}\label{indisc_set}
    \begin{enumerate}
	\item
	    Let $p$ be a stable type, $I$ an indiscernible sequence of
	    realizations of $p$. Then $I$ is an indiscernible set.
	\item
		Let $p$ be a stable type, $A$ a set, $I$ a sequence of
	    realizations of $p$ of length at least $(|A|+|T|)^+$. Then there
		exists an infinite subsequence $I'\subseteq I$, which is indiscernible over $A$.
    \end{enumerate}
\end{fct}

We now need to introduce the notion of a \emph{Morley sequence}. In
general, a \emph{Morley sequence} in a type $p \in \tS(A)$ is an
indiscernible sequence $I = \lseq{a}{i}{O}$ of realizations of $p$
such that $\tp(a_i/Aa_{<i})$ does not fork over $A$. Note that from
stationarity of stable types, Fact \ref{stationary}, it follows that
the only way to obtain a Morley sequence in a stable type $p$ over
an algebraically closed set $A$ is as follows: let $q$ be the unique
global nonforking extension of $p$. Define $\lseq{a}{i}{\om}$ such
that $a_i \models q \rest Aa_{<i}$. One still needs to make sure
that $I$ is indiscernible over $A$, but this comes for free:

\begin{fct}\label{fct:nonsplit_indisc}
    Let $q$ be a global type definable over a set $A = \acl(A)$.
    Define a sequence $I$ as described above. Then $I$ is
    indiscernible over $A$.
\end{fct}
\begin{prf}
    This is in fact true whenever $q$ is invariant under the action of
    $\Aut(\fC/A)$, see \cite{Sh:c}.
\end{prf}

\begin{dfn}
    Let $\lam$ be a cardinal. A \emph{block} of $\lam$ is a finite
    subset of $\lam$. For two blocks
    $u_1, u_2$ of
    $\lam$ we say that
    $u_1<u_2$ if $\max u_1 <
    \min u_2$.
\end{dfn}

\begin{prp}\label{prp:unique}(Strong Uniqueness)
    Let $p \in \tS_\Delta(A)$ be a minimal wide stable type, and let
    $I = \lseq{a}{\al}{\lam}$ be a Morley sequence in $p$. Then
    \begin{enumerate}
    \item
        $I$ is an indiscernible set over $A$.
    \item
        Let $u_i$ be mutually disjoint blocks of $\lam$ for $i<\om$ and $b_i
        \in \sum_{\al \in u_i} \setR a_\al$ with $\|b_i\| = 1$. Then $J =
        \lseq{b}{i}{\om}$ is an indiscernible set over $A$ and a
        Morley sequence in $p$.

        In particular, $\tp(J/A) = \tp(I/A)$.
    \end{enumerate}
\end{prp}
\begin{prf}
\begin{enumerate}
\item
    By stability (combine Fact \ref{fct:nonsplit_indisc} with Fact \ref{indisc_set}).
\item
    Let $p_\al = \tp(a_\al/Aa_{<\al})$. Fix $\al < \lam$. Note that
    $p_{\al+1}$ is a wide type extending $p_\al$. Let $B$ be a
    subspace of infinite dimension, isometric to $\ell_2$, whose unit sphere is contained in
    $p_{\al}^\fC$. We may assume that $a_\al \in \mathbf S(H)$. Clearly for
    all $a' \in \mathbf S(H)$ we have
    \begin{equation*}
    \setR a_\al + \setR a' \subseteq B
    \end{equation*}

    Moreover, if $r,r'\in \setR$ are such that $\|r+r'\|_2 = 1$, then
    for all $a' \in \mathbf S(H)$ we have
    $$ ra_\al + r'a' \in \mathbf S(H) $$

    Hence the following partial type over $Aa_{\le \al}$ is wide:

    $$ \pi(x) = \left\{p(ra_\al+r'x)\colon r,r' \in \setR, r^2+(r')^2 = 1 \right\}$$

    By Theorem \ref{thm:wideexist}, there exists a wide complete type  $p'(x)$ over $Aa_{\le\al}$ extending
    $\pi(x)$. Since $p'$
%
    clearly extends $p_\al$, by minimality we get $p' = p_{\al+1}$, so
    $a_\be \models p$ for all $\be > \al$. It follows by indiscernibility that for any
    $\be > \ga \ge \al$ and $r,r'$ with $\|r+r'\|_2 = 1$, we have
    $ra_\ga + r'a_\be \models p_\al$.




    Moreover, by clause (i), that is, since $I$ is an indiscernible set, it is easy to see that for any
    $\be > \ga \ge \al$ and $r,r'$ with $\|r+r'\|_2 = 1$, we have
    $ra_\ga + r'a_\be \models \tp(a_\al/Aa_{<\al}a_{>\be}) = \tp(a_\ga/Aa_{<\al}a_{>\be}) =
    \tp(a_\be/Aa_{<\al}a_{>\be})$. So denoting $a' = ra_\ga +
    r'a_\be$, we have that $I' = a_{<\al}{}^\frown a'^\frown a_{>\be}$
    is a Morley sequence in $p$.

    The case when $a'$ is a general block element is proven
    by induction. That is, suppose that $$a' = \sum_{i<n}r_ia_{\al_i} + r_na_{\al_n}$$
    such that $\sum_{i\le n} r_i^2 = 1$.
    By the induction hypothesis, denoting
    $$ r'' = \left\|\sum_{i<n}r_ia_{\al_i}\right\| = \sqrt{\sum_{i<n} r_i^2}$$ and
    $$ a'' = \frac{1}{r''}\sum_{i<n}r_ia_{\al_i} $$
    we have that the following sequence $$I'' = a_{<\al_0}{}^\frown
    a''{}^\frown a_{\al_n}{}^\frown a_{>\al_n}$$ is a Morley sequence in $p$. Note that
    $$(r'')^2 + r_n^2 = \sum_{i\le n}r_i = 1$$ so by the case
    $n = 2$ (which was our base case), the sequence $$ I' = a_{<\al_0}{}^\frown
    (r''a''+r_na_{\al_n}){}^\frown a_{>\al_n}$$ is a Morley
    sequence in $p$, as required.

    Now it is easy to deduce the general statement of Strong Uniqueness by
    induction on the number of blocks.
%
\end{enumerate}
\end{prf}


\begin{prp}\label{prp:l2}
    Let $p \in \tS_\Delta(A)$ be a minimal wide stable type, and let
    $I = \lseq{a}{\al}{\lam}$ be a Morley sequence in $p$. Then  $I$
    is isometric to the standard basis of $\ell_2$. In other words,
    for every $k<\om$ and $\lam_0, \ldots, \lam_{k-1} \in \setR$, we
    have
    $$\left\|\sum_{i<k}\lam_ia_i \right\|^2 =\sum_{i<k}|\lam_i|^2$$
\end{prp}
\begin{prf}
    Let $I'$  be (isometric to) the standard basis of an infinite dimensional
    $\ell_2$ space, $I' \subseteq p^\fC$. Since $I'$ can be chosen
    as large as we want, by stability there is $I \subseteq I'$
    indiscernible over $A$. Clearly $I$ is isometric to the standard
    basis of $\ell_2$. We need to show that $I$ is a Morley sequence
    over $A$. Let $H$ be the Hilbert space generated by $I$.

    Without loss of generality $A = \acl(A)$, so $p$ is stationary.
    Let $p^*$ be the global nonforking extension of $p$. Denote $I =
    \lseq{a}{i}{\om}$.

    Let $H_0$ be the subspace
    of $H$ generated by $Aa_0$. Note that all elements of the unit
    sphere of
    $(H_0)^\perp$ (the orthogonal complement in $H$), which is an
    infinite-dimensional Hilbert space, satisfy the partial type

    $$    \pi(x) = p(x)\bigcup \left\{\left\|\lam_0a_0+\lam x\right\|^2 =
    \lam_0^2+\lam^2\colon \lam_0,\lam \in
    \setR\right\}$$

    hence $\pi(x)$ is wide. By Theorem \ref{thm:wideexist},
    there exists $q\in \tS(Aa_{0})$ extending $\pi(x)$, which is wide.
    Since $q$ extends $p$ and is wide, by minimality of $p$ we
    have $q = p^*\rest Aa_0$. Let $b_0 = a_0, b_1 \models q$. Then
    $b_0,b_1$ start a Morley sequence in $p$, and as $q$ extends
    $\pi(x)$, we see that $\seq{b_0,b_1}$ is isometric to the
    standard basis of a two-dimensional Hilbert space.

    Now let $\lseq{b}{i}{\om}$ be a Morley sequence in $p$
    continuing $\lseq{b}{i}{2}$, and we show by induction on $n$
    that the sequence $\lseq{b}{i}{n}$ is isometric to the
    standard basis of an $n$-dimensional Hilbert space. Assuming
    that this holds for $\lseq{b}{i}{n}$, let us take care of
    $\lseq{b}{i}{n+1}$.

    Let $\lseq{\lam}{i}{n+1}$ be scalars in $\setR$. By the
    induction hypothesis we have

    \begin{equation}\tag{$\blacklozenge$}\label{equ:ind_hyp}
    \left\|\sum_{i<n}\lam_ib_i\right\|^2 = \sum_{i<n}\lam_i^2
    \end{equation}

    Denote $$\lam' = \sqrt{\sum_{i<n}\lam_i^2} \;\;\;\text{and}\;\; b'
    = \frac{1}{\lam'}\sum_{i<n}\lam_ib_i$$

    So $\|b'\| = 1$. By Strong Uniqueness (Proposition \ref{prp:unique}(ii)), the sequence $\seq{b',b_n}$ is
    a (2-element) Morley sequence in $p$. By the induction
    hypothesis again (or by the case $n=2$, which was our base
    case), we have
$$    \left\|\sum_{i<n+1}\lam_ib_i \right\|^2 = \left\|\lam'b' +
\lam_nb_n\right\|^2 = (\lam')^2+\lam_n^2 = \sum_{i<n}\lam_i^2 +
\lam_n^2 = \sum_{i<n+1}\lam_i^2$$

%
%
%

    which completes the induction step.

%
%
%
%
\end{prf}

\section{On Henson's Conjecture}\label{conclusions}

We recall that throughout this paper we are assuming that $K$ is an
elementary class of Banach spaces with extra-structure, $\fC$ its
monster model. In this section we will also assume that the language
of $K$ (which we denote by $L$) is \emph{countable}.

Let $M \in K$, $A \subseteq M$. We say that $M$ is \emph{prime} over
$A$ if whenever $A \subseteq N \in K$, there is an elementary
embedding $f \colon M \hookrightarrow N$ which is the identity on
$A$.

We now state some standard facts about non-separably categorical continuous theories \cite{ShUs837}, \cite{Ben:morley}. 

\begin{fct}\label{primes exist}
    Assume that $K$ is uncountably categorical. Let $A \subseteq
    \fC$. Then there exists a model $M \in K$ which is prime over
    $A$.
\end{fct}
\begin{prf}
    This is true in a more general context of $\aleph_0$-stable $K$. See section 4 of \cite{ShUs837} or
    \cite{Ben:morley}.
\end{prf}

\begin{fct}\label{stability}
    Assume that $K$ be uncountably categorical. Then $K$ is $\aleph_0$-stable, in particular stable.
\end{fct}

Recall that $K$ is stable if and only if every type in $\fC$ in the
language $L$ is stable.

\begin{fct}\label{Morley}(Morley's Theorem for continuous
logic, \cite{ShUs837, Ben:morley}).
    Assume that $K$ is uncountably categorical. Then $K$ is categorical in
    every uncountable density. Moreover, every non-separable model
    in $K$ is saturated.
\end{fct}

We are now ready to prove the main result of the paper.

\begin{thm}\label{thm:main}
    Let $K$ be an elementary class of Banach space with
    extra-structure, as defined in section 2, and assume that the
    language of $K$ is countable. Equivalently, assume that $T$ is a
    countable continuous theory whose monster model $\fC$ expands a
    Banach space $\fB$.

    Assume that $K$ (equivalently, $T$) is
    categorical in some
    uncountable density character. \underline{Then}: There is a separable model $M_0$ of $T$ and a wide type $p$ over $M_0$ such that
 	\begin{itemize}
	\item 
		Any Morley sequence in $p$ is isometric to the standard orthonormal basis of a Hilbert space;
	\item
		Any non-separable model of $T$ is prime over a Morley sequence in $p$. 
	\end{itemize}
	
	Specifically, if $M$ is a model of $T$ of uncountable density character $\lam$, then $M$ is prime over 
	a Morley sequence in $p$ of length $\lam$. 
	
	In particular, we have the following: Let $B_0$ be the Banach space that underlies $M_0$. Let $M\models T$ be of 
	uncountable density character $\lam$. Then there exists a spreading model $H$ of $M_0$ isometric to $\ell_2(\lam)$, 
	and $M$ is prime over $H$. 

%
\end{thm}
\begin{prf}
    By Fact \ref{primes exist}, let $\hat M_0$ be the prime model in $K$ (prime over $\emptyset$).
    By Theorem \ref{thm:wideexist}, there exists a wide type $\hat p_0$ over
    $\hat M_0$. By Fact \ref{stability}, $K$ is stable, in particular the type $\hat p_0$ is stable.
    By Corollary \ref{cor:minexist} (and e.g. Fact \ref{primes exist}, although it is not needed for this),
    there is a separable model $M_0 \in K$, $\hat M_0 \prec M_0$,
    and a minimal wide type extending $\hat p_0$.

    Now let $M \in K$ be of uncountable density $\lam$. By Fact
    \ref{Morley}, $M$ is $\lam$-saturated, so we may assume that $M_0
    \subseteq M$. By saturation again, there is a Morley sequence $I = \lseq{a}{i}{\lam}$ in
    $p_0$, $I \subseteq M$. Let $M'$ be a prime model over $I$ (Fact \ref{primes
    exist}). Then $M'$ has density $\lam$; since $K$ is categorical in $\lam$ by Fact \ref{Morley}, $M'$ is
    isometric to $M$.

    So $M$ is prime over a sequence isometric to a Morley sequence in
    $p_0$. The desired conclusion follows now from Proposition \ref{prp:l2}.
\end{prf}

We conclude with a few remarks and some possible directions for future research.

\begin{rmk}
	In Theorem \ref{thm:main}, one may assume that $M_0$ is the saturated separable model of $T$; however, not necessarily the prime model. 
	It would be interesting to find out whether an $\ell_2$ type exists over the prime model as well (we believe that the 
  answer ought to be positive).  
\end{rmk}

\begin{rmk}
  One could ask: in which sense have we shown that any $B \in K$ is determined by a Hilbert \emph{space}?
  It may seem from the way our main results are stated that in the more general case that $B$ is a Banach structure, that is, a continuous structure properly expanding a Banach space, what we have really proved is that its isometry type is controlled by an underlying Hilbert \emph{structure}; that is, a Hilbert space expanded with the additional structure on $B$. However, a close examination of Proposition \ref{prp:unique} indicates that the ``induced structure'' on $H$ is much simpler. In particular, the complete type (in the full language of $B$) of a singleton in $H$ is completely determined by its norm. It is of, course, problematic to talk about the induced structure on a non-definable set. Nevertheless, this suggests that the Hilbert space $H$ that ``controls'' $B$ is in a certain weak sense ``stably embedded''.
\end{rmk}

\begin{rmk}
  For a non-logician reader, interested in connections between stability (for example, as defined by Krivine-Maurey) and 
  definability of types (which is the notion that we have used here), we recommend an excellent short article by Ita\"{i} Ben Yaacov \cite{Ben:Gro}, where it is observed that definability of types in a stable theory (and moreover over a stable structure) can be derived from ``Grothndick's Criterion", characterizing sets of ``commuting'' functions (very much in the spirit of the Krivine-Maurey definition). 

  Since we work in a more general context of stable types, it does not exactly fall under Ben-Yaacov's treatment; nevertheless, \cite{Ben:Gro} can provide a good insight into how definability to the order property, especially for a functional analyst (alternative to the more classical model theoretic treatment in \cite{BU}).
\end{rmk}

\begin{rmk}\label{rmk:uniqueextension}
  The ``Hlibert'' type $p$ that we have found in Theorem \ref{thm:main} has several interesting properties, including
  a certain version of definability. Specifically, there exists an $M_0$-definable predicate $\theta(y)$ (which is this case 
  simply means a \emph{uniform} limit of functions of the form $y \longmapsto \|b+y\|$ with $b\in M_0$) such that \emph{for any} ultrafilter $\fU$ on $M_0$, any $M \in K$ containing $M_0$ and $p'=\Av_{qf}(\fU,M)$, we have that for any $a \in M$, $\|x+a\|^{p'} = \theta(a)$. 

  In other words, there exists a unique extension of $p$ to any elementary extension of $M_0$ which is given by an ultrafilter on $M_0$ (equivalently, finitely satisfiable in $M_0$). This is also the unique wide extension of $p$. 
  This unique extension is definable by an $M_0$-predicate.

  Note the strength of the latter definability assumption. From the general theory of spreading models, we know that for any ultrafilter $\fU$ on $M_0$, for any separable extension $M$ of $M_0$, denoting  $p'=\Av_{qf}(\fU,M)$, there is a sequence 
  $\lseq{b}{n}{\om}$ in $M_0$ such that for all $a\in M$, 

  $$\|x+a\|^{p'} = \lim_n\|b_n+a\|$$ 

  However:

  \begin{itemize}
    \item Even though $\|b_n+y\|$ are $M_0$-formulas, they do not necessarily converge uniformly;
    \item The choice of the sequence \lseq{b}{n}{\om} depends on the extension $M$;
    \item A ``good'' sequence as above is constructed by diagonalization, therefore generally exists only if $M$ is separable.
  \end{itemize}

  By contrast, the definition $\theta(y)$ is a uniform limit of $M_0$-formulas, and it works uniformly for any extension of $M_0$, separable or otherwise. 

  One can see the above uniformity in several more explicit ways. First, by an easy application of Mazur's Lemma (see
  \cite{Ben:Gro}, Corollary 7) one can assume that $\theta(y)$ is given by a uniform limit of averages of formulas of the form 
$\|b_n+y\|$, where $b_n$ is a certain approximating sequence for $p$.

  In addition, if one is willing to give up the fact that it is defined \emph{over $M_0$}, one can obtain a very 
  explicit formula for it. Specifically, let \lseq{b}{n}{\om} be any Morley sequence in $p$. Then for any ultrafilter $\fU$ on $M_0$, for any elementary extension $M$ of $M_0$, denoting  $p'=\Av_{qf}(\fU,M)$, 
  for all $a\in M$, we have $\|x+a\|^{p'} = \lim\|b_n+a\|$. Moreover, the uniformity can be expressed in the following way: for any $\eps>0$ there exists $N<\om$ such that 

  $$\|x+a\|^{p'} \sim_\eps \max_{W\subset 2N, |W| = N+1}\min_{n \in W} \|b_n + a\|$$

  All this follows from the continuous analogue of classical stability theory; see \cite{BU}.

  It would still be interesting to know whether yet stronger versions of definability are true. In particular, we ask: is the Hilbert space $H$, given as the spreading model of $p$, a type-definable \emph{set}? Even more, is $H$ a zero set of a definable predicate?

  In an earlier version of this article we asked whether $H$ can be assumed to be \emph{definable}; recently C. Ward Henson has described to us a large class of examples (constructed by Raynaud and himself) that show that in general this is too much to hope for. 
\end{rmk}

\appendix

\section{More on existence of minimal types}\label{app:a}
In this appendix, we give an alternative proof of the fact that if $\fC$ is q.f. stable, then any wide type $\pi$ over a model 
$M$ can be extended to a minimal wide type $\pi'$ over $M \cup A$ where $A$ is countable. This proof gives less information that the one in Section 4 (for example, it does not immediately imply that a sequence in the 
global wide extension of $\pi'$ gives rise to a spreading model of $\Span(M\cup A)$), this is why we decided to leave it for the appendix. However, it also has an important advantage: it
does not invoke model theoretic ``black boxes'' such as definability of stable types (equivalently, ``Grothendick's criterion''), thereby making it  more illuminating and transparent for a non-expert. We believe that it has practical 
relevance even for a model theorist; for example, it may help addressing the question of definability, stated at the end of the previous section. 
This approach was inspired by a talk given by Angus Macintyre in the Logic Seminar at Centro de Matem\'{a}tica e 
Aplica\c \~{o}es Fundamentais in Lisbon.  

Throughout this section, we work in the simple context outlined in section 2.1. That is, we assume that $K$ is an elementary class of Banach spaces (with no extra-structure). In addition, we assume that its monster model $\fC$ is quantifier free stable (equivalently, every $M$ in $K$ is Krivine-Maurey stable -- this property is called ``super-stability'' by Banach space theories). The context can be
significantly generalized (and the assumptions weakened), but since the appendix is mostly written for non-logicians (and 
provides an alternative proof of a result that appears in full generality in the main body of the article), we do not see much 
point in doing so. Hence the term ``types'' will refer exclusively to quantifier free types, and ``stability'' will mean 
Krivine-Maurey stability. However, the reader could easily generalize the proofs to the more general context of the article, should they be interested in doing so. 

The proof  relies on the following important characterization of stability. The proof is not hard, and can be found in various sources
mentioned in the article; we do not include it here. In fact, this is the characterization of stability that follows most directly from non-separable categoricity, and is proven in \cite{Ben:morley} and \cite{ShUs837}.

First, a definition: given $A \subseteq \fC$ and two (partial) types $p,q$ over $M$, 
we define $d(p,q)$, the distance between $p$ and $q$, to be the infimum of $\|a-b\|$ where $a,b$ range over realizations (in $\fC$) of $p$ and $q$, respectively). 

\begin{fct}
  $\fC$ is stable if and only if given any separable subspace $M$, there is a separable $N \in K$ containing $M$ such that 
  all types over $M$ are realized in $N$.
\end{fct}

In other words, the space of types over a separable (equivalently, countable) susbset of $\fC$ is separable (with respect to the metric $d$ defined above). 

Note that if $\lseq{p}{i}{\lam}$ is a sequence of types over $M$ such that for some $\eps>0$ for any $i<j$ there is $a\in M$ such that 

$$\left| \|x+a\|^{p_i} - \|x+a\|^{p_j}\right| \ge \eps$$ 

then the density character of $S_{qf}(M)$ with respect to $d$ is at least $\lam$. This is because if $b_i$ is a realization of $p_i$, then 
$$ \|b_i-b_j\| \ge \left| \|b_i+a\| - \|b_j+a\| \right| \ge \eps$$

In particular, 

\begin{fct}\label{fct:instability}
Let $M$ be separable. Assume that 
there are $\eps>0$ and  uncountably many types $p_i$ over $M$ so that or any $i<j$ there is $a\in M$ such that 

$$\left| \|x+a\|^{p_i} - \|x+a\|^{p_j}\right| \ge \eps$$ 

Then $\fC$ is unstable. 
\end{fct}

In order to simplify the proof of the main theorem of this section, we will make the following definition:

\begin{dfn}
  Let $\pi(x)$ be a wide type, and $\eps>0$. We say that $\pi$ is \emph{$\eps$-explicitly non-minimal} if there exists $a \in \fC$
  and $r<s \in \setR$ with $|s-r|\ge \eps$ such that both 
  $\pi_{\seq{0}}=\pi(x)\cup\set{\|x-a\|\le r}$ and 
  $\pi_{\seq{1}} = \pi(x)\cup\set{\|x-a\|\ge s}$ are wide types. 

  We call the negation of the above notion simply \emph{$\eps$-minimal}.
 \end{dfn}

\begin{obs}\label{obs:explicitminimal}
Clearly, if $\pi(x)$ is a wide type which is not wide minimal, it is $\eps$-explicitly non-minimal for some $\eps>0$.

In other words, if $\pi$ is $\eps$-minimal for all $\eps>0$, then it is wide minimal.
\end{obs}

We also introduce the following notation: let $\pi_0, \pi_1$ be types such that for some $a \in \fC$ and $r<s \in \setR$ with $|s-r|\ge \eps$ we have $\set{\|x-a\|\le r} \subseteq \pi_0$ and 
  $\set{\|x-a\|\ge s}\subseteq \pi_1$. Then we say that the $\pi_0$ and $\pi_1$ are explicitly of distance at least $\eps$ from each other, and write $D(\pi_0,\pi_1)\ge\eps$. 

  Using this notation, we can restate Fact \ref{fct:instability} in the following convenient form:

  \begin{obs}\label{obs:instability}
    Let $M$ be separable. Assume that there are $\eps>0$ and  uncountably many types $p_i$ over $M$ so that or any $i<j$ we have $D(p_i,p_j)\ge\eps$. Then $\fC$ is unstable.
  \end{obs}

\begin{rmk}
  One can state a ``local'' version of the above Observation; that is, if a type $p$ has uncountably many such extensions, then $p$ is unstable. This provides an alternative root towards the proof of Corollary \ref{cor:minexist}.
\end{rmk}

\begin{lem}\label{lem:epsminimal}
  Assume that $\fC$ is stable, $\pi(x)$ a wide type (or just that $\pi$ is a stable wide type), and let $\eps>0$. Then there exists a finite set $B$ and an extension $\pi_\eps(x)$ of $\pi$ to a wide type over $B$ such that $\pi_\eps$ is $\eps$-minimal.
\end{lem}
\begin{prf}
  If $\pi$ is $\eps$-minimal, we are done. Otherwise, for some $a\in M$ and $r<s \in \setR$, 
  where $|s-r|\ge\eps$,both $\pi_{\seq{0}}=\pi(x)\cup\set{\|x-a\|\le r}$ and 
  $\pi_{\seq{1}} = \pi(x)\cup\set{\|x-a\|\ge s}$ are wide types. Clearly $D(\pi_{\seq{0}},\pi_{\seq{1}})\ge\eps$. 
  If $\pi_{\seq{0}}$ is $\eps$-explicitly non-minimal, we can construct extensions $\pi_{\seq{00}}$ and $\pi_{\seq{01}}$ of
   $\pi_{\seq{0}}$ with $D(\pi_{\seq{00}},\pi_{\seq{01}})\ge\eps$. Similarly, if $\pi_{\seq{1}}$ is $\eps$-explicitly non-minimal, we can construct extensions $\pi_{\seq{10}}$ and $\pi_{\seq{11}}$ of $\pi_{\seq{1}}$ with 
   $D(\pi_{\seq{10}},\pi_{\seq{11}})\ge\eps$.

   This is where stability comes in. If this construction could continue, we would get types $\pi_\eta$ for $\eta \in {}^{\om>}2$ (where ${}^{\om>}2$ denotes 
   the binary tree) 
   such that if $\eta \triangleleft \nu \implies \pi_\eta \subseteq \pi_\nu$, and if $\eta,\nu$ are incomparable, then 
   $D(\pi_\eta,\pi_\nu)\ge\eps$. Taking unions of the types along the branches, we obtain continuum many types over a countable set of parameters (since the binary tree has countably many nodes, and at each split we only added
   one new element), which are 
   pairwise explicitly of distance $\ge \eps$. This contradicts Observation \ref{obs:instability}.

   Since the construction fails, there is $\eta \in {}^{\om>} 2$ which is $\eps$-minimal. Letting $B$ be the collection of elements that were added as parameters along this (finite) branch, we are done. 

\end{prf}

As mentioned above, in the following theorem we assume for simplicity of exposition that $\fC$ is stable; however, a suitable local analogue is easy to state (and prove). 

\begin{thm}($\fC$ stable)\label{thm:appendix}
  Let $M \in K$, and let $\pi(x)$ be a partial wide type over $M$. Then there exists a countable set $A$ such that there
  exists a minimal wide type $\pi'(x)$ over $M\cup A$ extending $\pi$. 
\end{thm}

\begin{prf}
  Assume $\pi(x)$ is not minimal. By Observation \ref{obs:explicitminimal}, it $\eps$-explicitly non-minimal for some $\eps>0$.
  
  Denote $\delta_n = \frac{\eps}{2^n}$. Construct by induction on $n$ a type $\pi_n(x)$ such that:

  \begin{itemize}
    \item $n<m \implies \pi \subseteq \pi_n \subseteq \pi_m$
    \item  $n<m \implies \dom(\pi_m)\setminus\dom(\pi_n)$ is finite. 
    \item $\pi_n$ is wide and $\delta_n$-minimal
  \end{itemize}

  This can be easily accomplished by repeatedly applying Lemma \ref{lem:epsminimal}.

  Now the union $\pi'(x) = \bigcup_{n<\om}\pi_n(x)$ is clearly as required.

\end{prf}

For the sake of completeness of presentation we now outline a path to the proof of the main theorem of this paper that does not go through the ``black box'' of definable types. 

Let $K$ be a non-separably categorical class of Banach spaces. Let $\hat{M_0}$ be the prime model of $K$, and let $\pi$ be a wide type over $\hat{M_0}$. $K$ is stable, hence by the above Theorem, there is countable extension of $\hat{M_1}$ of $\hat{M_0}$ over which $\pi$ has a wide minimal extension. In fact, by Proposition \ref{prp:wideextend} we may assume that this extension is a complete type over $\hat{M_1}$; call it $\hat{p}$. Take a countable Morley sequence $I=\lseq{a}{n}{\om}$ in $\hat{p}$; we have seen in Proposition \ref{prp:l2} that it is isometric to the standard basis of $\ell_2$. Take a separable model $M_0 \in K$ containing $\hat{M_1}$ and the sequence $I$ (one could take the prime model over $\hat{M_1}\cup I$; alternatively, the separable saturated model of $K$ would work). Let $p$ be the unique wide extension of $\hat{p}$ to a complete type over $M_0$. It is now not too hard to see that $p$ satisfies all the requirements of the type whose existence is postulated in Theorem \ref{main_theorem}. In particular, its unique wide extension to any super-structure of $M_0$ is finitely satisfiable in $M_0$ (and indeed in the sequence $I$), and in fact the sequence $I$ approximates any spreading model generated by $p$. Moreover, the approximation is uniform in the sense explained in Remark \ref{rmk:uniqueextension}, which leads to definability of the type $p$ over $M_0$ in the way hinted at there (some more local stability theory, as developed in \cite{BU} is required in order to understand why this is the case). 

The proof outlined here still gives less information than the approach in the main body of this article. In particular, the proof of Theorem \ref{thm:main} yields definability over a smaller structure (which may play a role in the attempt to obtain a minimal type over the prime model), and suggests various other generalizations that will be explored elsewhere. However, we believe that the proof of  Theorem \ref{thm:appendix} as presented here (more precisely, the proof of Lemma \ref{lem:epsminimal}) gives a certain insight into the importance and the meaning of stability in this context, and therefore was worth including.

\bibliography{common.bib}

\def\cprime{$'$}
\begin{thebibliography}{BYBHU08}

\bibitem[BL71]{BaLa}
J.~T. Baldwin and A.~H. Lachlan.
\newblock On strongly minimal sets.
\newblock {\em J. Symbolic Logic}, 36:79--96, 1971.

\bibitem[BL00]{BL}
Yoav Benyamini and Joram Lindenstrauss.
\newblock {\em Geometric nonlinear functional analysis. {V}ol. 1}, volume~48 of
  {\em American Mathematical Society Colloquium Publications}.
\newblock American Mathematical Society, Providence, RI, 2000.

\bibitem[BS74]{BrSu}
Antoine Brunel and Louis Sucheston.
\newblock On {$B$}-convex {B}anach spaces.
\newblock {\em Math. Systems Theory}, 7(4):294--299, 1974.

\bibitem[BY05]{Ben:morley}
Itay Ben~Yaacov.
\newblock Uncountable dense categoricity in cats.
\newblock {\em J. Symbolic Logic}, 70(3):829--860, 2005.

\bibitem[BY14]{Ben:Gro}
Ita{\"{\i}} Ben~Yaacov.
\newblock Model theoretic stability and definability of types, after {A}.
  {G}rothendieck.
\newblock {\em Bull. Symb. Log.}, 20(4):491--496, 2014.

\bibitem[BYBHU08]{BBHU}
Ita\"{i} Ben~Yaacov, Alexander Berenstein, C.~Ward Henson, and Alexander
  Usvyatsov.
\newblock Model theory for metric structures.
\newblock In {\em Model Theory with Applications to Algebra and Analysis (II)},
  volume 350 of {\em Lecture Notes series of the London Mathematical Society},
  pages 315--428. Cambridge University Press, Cambridge, 2008.

\bibitem[BYU10]{BU}
Ita{\"{\i}} Ben~Yaacov and Alexander Usvyatsov.
\newblock Continuous first order logic and local stability.
\newblock {\em Trans. Amer. Math. Soc.}, 362(10):5213--5259, 2010.

\bibitem[Dvo61]{Dv}
Aryeh Dvoretzky.
\newblock Some results on convex bodies and {B}anach spaces.
\newblock In {\em Proc. {I}nternat. {S}ympos. {L}inear {S}paces ({J}erusalem,
  1960)}, pages 123--160. Jerusalem Academic Press, Jerusalem, 1961.

\bibitem[GOU13]{GOU}
Dar\'{i}o Garc\'{i}a, Alf Onshuus, and Alexander Usvyatsov.
\newblock Generic stability, forking, and thorn-forking.
\newblock {\em Trans. Amer. Math. Soc.}, 365(1):1--22, 2013.

\bibitem[Hei80]{Hei}
Stefan Heinrich.
\newblock Ultraproducts in {B}anach space theory.
\newblock {\em J. Reine Angew. Math.}, 313:72--104, 1980.

\bibitem[Hen76]{Hen}
C.~Ward Henson.
\newblock Nonstandard hulls of {B}anach spaces.
\newblock {\em Israel J. Math.}, 25(1-2):108--144, 1976.

\bibitem[HI02]{HI}
C.~Ward Henson and Jos{\'e} Iovino.
\newblock Ultraproducts in analysis.
\newblock In {\em Analysis and logic (Mons, 1997)}, volume 262 of {\em London
  Math. Soc. Lecture Note Ser.}, pages 1--110. Cambridge Univ. Press,
  Cambridge, 2002.

\bibitem[HP11]{HP}
Ehud Hrushovski and Anand Pillay.
\newblock On {NIP} and invariant measures.
\newblock {\em J. Eur. Math. Soc. (JEMS)}, 13(4):1005--1061, 2011.

\bibitem[HR16]{HeRa}
C.~Ward Henson and Yves Raynaud.
\newblock Asymptotically {H}ilbertian modular {B}anach spaces: examples of
  uncountable categoricity.
\newblock {\em Comment. Math.}, 56(1):119--144, 2016.

\bibitem[Iov99a]{Iov:stabI}
Jos{\'e} Iovino.
\newblock Stable {B}anach spaces and {B}anach space structures. {I}.
  {F}undamentals.
\newblock In {\em Models, algebras, and proofs ({B}ogot\'a, 1995)}, volume 203
  of {\em Lecture Notes in Pure and Appl. Math.}, pages 77--95. Dekker, New
  York, 1999.

\bibitem[Iov99b]{Iov:stabII}
Jos{\'e} Iovino.
\newblock Stable {B}anach spaces and {B}anach space structures. {II}. {F}orking
  and compact topologies.
\newblock In {\em Models, algebras, and proofs ({B}ogot\'a, 1995)}, volume 203
  of {\em Lecture Notes in Pure and Appl. Math.}, pages 97--117. Dekker, New
  York, 1999.

\bibitem[Iov05]{Iov:def}
Jos{\'e} Iovino.
\newblock Definable types over {B}anach spaces.
\newblock {\em Notre Dame J. Formal Logic}, 46(1):19--50 (electronic), 2005.

\bibitem[KM81]{KM}
J.-L. Krivine and B.~Maurey.
\newblock Espaces de {B}anach stables.
\newblock {\em Israel J. Math.}, 39(4):273--295, 1981.

\bibitem[LP79]{LP}
Daniel Lascar and Bruno Poizat.
\newblock An introduction to forking.
\newblock {\em J. Symbolic Logic}, 44(3):330--350, 1979.

\bibitem[Mac71]{Mac}
Angus Macintyre.
\newblock On {$\omega \sb{1}$}-categorical theories of fields.
\newblock {\em Fund. Math.}, 71(1):1--25. (errata insert), 1971.

\bibitem[Mil92]{Mil}
V.~Milman.
\newblock Dvoretzky's theorem---thirty years later.
\newblock {\em Geom. Funct. Anal.}, 2(4):455--479, 1992.

\bibitem[Mor65]{Mor}
Michael Morley.
\newblock Categoricity in power.
\newblock {\em Trans. Amer. Math. Soc.}, 114:514--538, 1965.

\bibitem[Pes06]{Pes}
Vladimir Pestov.
\newblock {\em Dynamics of infinite-dimensional groups}, volume~40 of {\em
  University Lecture Series}.
\newblock American Mathematical Society, Providence, RI, 2006.
\newblock The Ramsey-Dvoretzky-Milman phenomenon, Revised edition of {{\i}t
  Dynamics of infinite-dimensional groups and Ramsey-type phenomena} [Inst.
  Mat. Pura. Apl. (IMPA), Rio de Janeiro, 2005; MR2164572].

\bibitem[PT11]{PT}
Anand Pillay and Predrag Tanovi{\'c}.
\newblock Generic stability, regularity, and quasiminimality.
\newblock In {\em Models, logics, and higher-dimensional categories}, volume~53
  of {\em CRM Proc. Lecture Notes}, pages 189--211. Amer. Math. Soc.,
  Providence, RI, 2011.

\bibitem[She90]{Sh:c}
S.~Shelah.
\newblock {\em Classification theory and the number of nonisomorphic models},
  volume~92 of {\em Studies in Logic and the Foundations of Mathematics}.
\newblock North-Holland Publishing Co., Amsterdam, second edition, 1990.

\bibitem[She04]{Sh715}
Saharon Shelah.
\newblock Classification theory for elementary classes with the dependence
  property---a modest beginning.
\newblock {\em Sci. Math. Jpn.}, 59(2):265--316, 2004.
\newblock Special issue on set theory and algebraic model theory.

\bibitem[SU11]{ShUs837}
Saharon Shelah and Alexander Usvyatsov.
\newblock Model theoretic stability and categoricity for complete metric
  spaces.
\newblock {\em Israel J. Math.}, 182:157--198, 2011.

\bibitem[Usv09]{Us}
Alexander Usvyatsov.
\newblock Generically stable types in dependent theories.
\newblock {\em J. Symbolic Logic}, 74(1):216--250, 2009.

\bibitem[Zil93]{Zil}
Boris Zilber.
\newblock {\em Uncountably categorical theories}, volume 117 of {\em
  Translations of Mathematical Monographs}.
\newblock American Mathematical Society, Providence, RI, 1993.
\newblock Translated from the Russian by D. Louvish.

\end{thebibliography}
\bibliographystyle{alpha}
\end{document}